\documentclass[12pt,a4]{amcjou}
\usepackage{amssymb}
\usepackage{latexsym}
\usepackage{color}
\begin{document}

\begin{frontmatter}   

\titledata{Polytopes with groups of type $PGL_2(q)$}           
{}                 

\authordata{Dimitri Leemans}            
{Universit\'e Libre de Bruxelles,
D\'epartement de Math\'ematiques - C.P.216,
Boulevard du Triomphe,
B-1050 Bruxelles}    
{dleemans@ulb.ac.be}                     
{Supported by the``Communaut\'e Fran\c caise de Belgique - Actions de Recherche Concert\'ees"}          

\authordata{Egon Schulte}            
{Northeastern University,
Department of Mathematics,
360 Huntington Avenue,
Boston, MA 02115, USA} 
{schulte@neu.edu}
{Supported by NSA-grant H98230-07-1-0005}                                       


\begin{abstract}
There exists just one regular polytope of rank larger than 3 whose full automorphism group is a projective general linear group $PGL_2(q)$, for some prime-power $q$. 
This polytope is the $4$-simplex and the corresponding group is $PGL_2(5)\cong S_5$.
\end{abstract}

\end{frontmatter}   

\section{Introduction}
\label{intro}

In an earlier paper~\cite{LS2005a} we determined the projective special linear groups $L_2(q)$ (often denoted $PSL_2(q)$), for $q$ a prime power, which occur as full automorphism groups of abstract regular polytopes of rank $4$ or higher.  We found that the only groups possible are $L_2(11)$ and $L_2(19)$, associated with two locally projective regular polytopes of rank $4$, namely Gr\"unbaum's 11-cell of type $\{3,5,3\}$ and Coxeter's 57-cell of type $\{5,3,5\}$, respectively (see \cite{Cox1982,Cox1984,Gr1977}). 

In the present paper we investigate the projective general linear groups $PGL_2(q)$, for $q$ a prime power. As for $L_2(q)$, these groups cannot be the full automorphism groups of abstract regular polytopes of rank $5$ or higher. Moreover, we shall see that if $PGL_2(q)$ is the full automorphism group of an abstract regular polytope of rank $4$, then necessarily $q=5$ and the polytope is the $4$-simplex with group $PGL_2(5)\cong S_5$.  This is in stark contrast to the situation in rank $3$, a phenomenon already observed for the groups $L_2(q)$. The exceptional case is based on the sporadic isomorphism between $PGL_2(5)$ and $S_5$. There is a wealth of regular polyhedra or maps on surfaces with groups (or rotation subgroups) isomorphic to $PGL_2(q)$ or $L_2(q)$ (see, for example,  McMullen, Monson \& Weiss~\cite{MMW1993}, Glover \& Sjerve~\cite{GS1985}, Conder~\cite{Co1990}).  In particular it is known that $PGL_2(q)$ is the automorphism group of a regular polyhedron for any $q>2$ (see Sjerve \& Cherkassoff~\cite{SC94}); the same is true for $L_2(q)$ when $q \neq 2, 3, 7$ or~$9$. 

\section{Basic notions}
\label{def}

For general background on (abstract) regular polytopes and C-groups we refer to McMullen \& Schulte~\cite{MS2002}. 

A {\it polytope} $\cal P$ is a ranked partially ordered set whose elements are called {\it faces}. A polytope $\cal P$ of rank $n$ has faces of ranks $-1,0,\ldots,n$; the faces of ranks $0$, $1$ or $n-1$ are also called {\it vertices}, {\it edges} or {\it facets}, respectively. In particular, $\cal P$ has a smallest and a largest face, of rank $-1$ and $n$, respectively. 
A {\it flag} of $\cal P$ is a maximal subset of pairwise incident faces of $\cal P$.
Each flag of $\cal P$ contains $n+2$ faces, one for each rank.  In addition to being locally and globally connected (in a well-defined sense), $\cal P$ is {\it thin\/}; that is, for every flag and every $j = 0,\ldots,n-1$, there is precisely {\it one} other ({\em $j$-adjacent\/}) flag with the same faces except the $j$-face. A polytope of rank $3$ is a {\em polyhedron\/}.  A polytope $\cal P$ is {\it regular\/} if its (combinatorial automorphism) group $\Gamma ({\cal P})$ is transitive on the flags. If $\Gamma({\cal P})$ has exactly two orbits on the flags such that adjacent flags are in distinct orbits, then $\cal P$ is said to be {\it chiral\/}.  

The groups of regular polytopes are string C-groups, and vice versa.  A {\em C-group\/} of {\em rank\/} $n$ is a group $G$ generated by pairwise distinct involutions $\rho_0,\ldots,\rho_{n-1}$ which satisfy the following {\it intersection property\/}:
\[  \langle \rho_j \mid j \in J\rangle \cap \langle \rho_j \mid j \in K\rangle = 
\langle\rho_j \mid j \in J\cap K\rangle \quad \quad
(\mathrm{for~} J, K \subseteq \{0,\ldots,n-1\}).\]
Moreover, $G$, or rather $(G, \{\rho_0,\ldots, \rho_{n-1}\})$, is a {\em string\/} C-group if the underlying Coxeter diagram is a string diagram; that is, if the generators satisfy the relations
\[ (\rho_j\rho_k)^2 = 1 \quad \quad
(\mathrm{for~} 0 \leq j < k -1 \leq n-2).\]
Each string C-group $G$ determines (uniquely) a regular $n$-polytope $\cal P$ with automorphism group $G$. The {\em $i$-faces} of $\cal P$ are the right cosets of the distinguished subgroup $G_i := \langle \rho_{j} \mid j \neq i \rangle$ for each $i = 0,1,\ldots, n-1$, and two faces are incident just when they intersect as cosets;   formally,  we must adjoin two copies of $G$ itself, as the (unique) $(-1)$- and $n$-faces of $\cal P$.  Conversely, the group $\Gamma({\cal P})$ of a regular $n$-polytope $\cal P$ is a string C-group, whose generators $\rho_j$ map a fixed, or {\em base\/}, flag $\Phi$ of $\cal P$ to the $j$-adjacent flag $\Phi^{j}$ (differing from $\Phi$ in the $j$-face).

We denote by $D_{2n}$ a dihedral group of order $2n$ and by $Z_n$ a cyclic group of order $n$.

\section{The subgroup structure of $PGL_2(q)$}
\label{pslsub}

Our first proposition reviews the subgroup structure of $L_2(q)$ and may be found in Dickson~\cite{Dic58} (but was first obtained in papers by Moore~\cite{Moo1904} and Wiman~\cite{Wi}); see also Huppert~\cite[Ch. II, \S 8]{Hup67} for a weaker version, and Kantor~\cite{ka} for interesting historical information about these groups.

\begin{prop}\label{sub}

The group $L_2(q)$ of order $\frac{q  (q^2 - 1)}{(2,q-1)}$, where $q=p^r$ with $p$ a
prime and $r$ a positive integer, contains only the following subgroups:

\begin{enumerate}
\item Elementary abelian subgroups of order $q$, denoted by $E_q$.

\item Cyclic subgroups $Z_d$, for all divisors $d$ of $\frac{(q\pm 1)}{(2,q-1)}$.

\item $\frac{q(q^2 -1)}{2d (2,q-1)}$ dihedral subgroups groups $D_{2d}$,
for all divisors $d$ of $\frac{(q\pm 1)}{(2,q-1)}$ with $d>2$. The number of conjugacy
classes of these subgroups is $1$ if $\frac{(q\pm 1)}{d (2,q-1)}$ is odd, and $2$ if
it is even.

\item For $q$ odd, $\frac{q(q^2 -1)}{12 (2,q-1)}$ dihedral subgroups $D_4$ (Klein 4-groups). The number of conjugacy classes of these subgroups is $1$ if  $q \equiv \pm 3 (8)$ and $2$ if $q \equiv \pm 1 (8)$. For $q$ even, the subgroups $D_4$ are listed under family~(5).

\item Elementary abelian subgroups of order $p^s$, denoted by $E_{p^s}$, for all natural number $s$ such that 
$1\leq s \leq r-1$.

\item Subgroups $E_{p^s}\!:\!Z_h$, each a semidirect product of an elementary abelian subgroup $E_{p^s}$ and a cyclic subgroup $Z_h$, for all natural numbers $s$ such that $1\leq s \leq r$ and all divisors $h$ of $\frac {p^k -1}{(2,1,1)}$, where $k=(r,s)$ and $(2,1,1)$ is defined as 2, 1 or 1 according as $p>2$ and $\frac{r}{k}$ is even, $p>2$ and $\frac{r}{k}$ is odd, or $p=2$.

\item For $q$ odd or $q=4^m$,\  alternating groups $A_4$, of order 12.

\item For $q \equiv \pm 1 (8)$, symmetric groups $S_4$, of order 24.

\item For $q \equiv \pm 1 (5)$ or $q=4^m$, alternating groups $A_5$, of order 60. For $q \equiv 0 (5)$, the groups $A_5$ are listed under family~(10).

\item $\frac{q(q^2 -1)}{ p^w (p^{2w}-1)}$ groups $L_2(p^w)$, for all divisors $w$ of $r$. The number of conjugacy classes of these subgroups is $2$, $1$ or $1$ according as $p>2$ and $\frac {r}{w}$ is even, $p>2$ and $\frac {r}{w}$ is odd, or $p=2$.

\item Groups $PGL_2(p^w)$, for all $w$ such that $2w$ is a divisor of $r$.
\end{enumerate}
\end{prop}

Observe here that when $q$ is even, $PGL_{2}(q) \cong L_{2}(q)$, so that family~(11) of Proposition~\ref{sub} really is a subfamily of family~(10) in that case.

We frequently require properties of the subgroup lattice of $PGL_2(q)$.  Two key observations employed in our proofs are that $PGL_2(q)$ can be viewed as a subgroup of $L_2(q^2)$ and that $PGL_2(q)$ has a unique subgroup isomorphic to $L_2(q)$.  This allows to extract the list of subgroups of $PGL_2(q)$ from the list of subgroups of $L_2(q^2)$, leading to the following result taken from Cameron, Omidi \& Tayfeh-Rezaie~\cite[\S 3]{COT2006} (see also Moore~\cite{Moo1904}). If $q=p^r$, we let $\epsilon = \pm 1$ be determined by $q \equiv \epsilon$ 
mod~$4$.

\begin{prop}\label{subPGL}

The subgroups of $PGL_2(q)$, where $q=p^r$ with $p$ an odd
prime and $r$ a positive integer, are as follows.

\begin{enumerate}
\item Two conjugacy classes of cyclic subgroups $Z_2$:\  one class consisting of $q(q+\epsilon)/2$ subgroups $Z_2$ contained in the subgroup $L_2(q)$, the other of $q(q-\epsilon)/2$ subgroups $Z_2$ not contained in $L_2(q)$.
\item One conjugacy class of $q(q\mp \epsilon)/2$ cyclic subgroups $Z_d$, where $d\mid q\pm \epsilon$ and $d>2$.
\item Two conjugacy classes of dihedral subgroups $D_4$:\ one class consisting of $q(q^2-1)/24$ subgroups $D_4$ contained in the subgroup $L_2(q)$, the other of $q(q^2-1)/8$ subgroups $D_4$ not contained in $L_2(q)$.
\item Two conjugacy classes of dihedral subgroups $D_{2d}$, where $d\mid\frac{q\pm\epsilon}{2}$ and $d>2$:\ one class consisting of $q(q^2-1)/(4d)$ subgroups $D_{2d}$ contained in the subgroup $L_2(q)$, the other of $q(q^2-1)/(4d)$ subgroups $D_{2d}$ not contained in $L_2(q)$.
\item One conjugacy class of $q(q^2-1)/(2d)$ dihedral subgroups $D_{2d}$, where $(q\pm\epsilon)/d$ is an odd integer and $d>2$.
\item $q(q^2-1)/24$ subgroups isomorphic to $A_4$, $q(q^2-1)/24$ subgroups isomorphic to $S_4$, and $q(q^2-1)/60$ subgroups isomorphic to $A_5$ when $q=\pm 1(10)$. There is only one conjugacy class of each of these types of subgroups, and all subgroups lie in the subgroup $L_{2}(q)$ except for those of type $S_4$ when $q=\pm 3 (8)$.
\item One conjugacy class of $p^r(p^{2r}-1)/(p^m(p^{2m}-1))$ subgroups $L_{2}(p^m)$, where $m\mid r$.
\item The subgroups $PGL(2,p^m)$ where $m\mid r$.
\item A semidirect product of the elementary abelian group of order $p^m$ with $m\leq r$ and the cyclic group of order $d$ with $d\mid q-1$ and $d\mid p^m-1$.
\end{enumerate}
\end{prop}

Note in particular that $L_2(q)$ and $PGL_2(q)$ do have dihedral subgroups $D_{2p}$, even though they are not explicitly listed among the dihedral groups in items (3,4) of Proposition~\ref{sub} or items (3,4,5) of Proposition~\ref{subPGL}, respectively; however, they do show up among the groups in items (6) or (9) of the respective Propositions.

On the other hand, $L_2(q)$ and $PGL_2(q)$ do not have subgroups isomorphic to $D_{4p}$. 
This makes it impossible for the subgroup $G_{23}=\langle\rho_{0},\rho_{1}\rangle$ (respectively $G_{01}=\langle\rho_{2},\rho_{3}\rangle$) in the next sections to be isomorphic to $D_{2p}$ because otherwise, $G_{2}=\langle\rho_{0},\rho_{1},\rho_{3}\rangle$ (respectively $G_{1}=\langle\rho_{0},\rho_{2},\rho_{3}\rangle$) would have to be isomorphic to $D_{4p}$.
Observe that this is what permits us to consider only dihedral subgroups of the maximal dihedral subgroups in the proof of Lemma 7 of~\cite{LS2005a} and should have been pointed out in \cite{LS2005a}. A similar comment applies to Lemma~\ref{lemma4}.

\section{$PGL_2(q)$ acting on polytopes}
\label{pslacts}

In this section, we assume that $G$ is a group isomorphic to $PGL_2(q)$, with $q=p^r$. The prime $p=2$ is special: in fact, if $p=2$ then $PGL_2(q) \cong L_2(q)$ and we may simply appeal to our paper~\cite{LS2005a} about polytopes with groups $L_2(q)$ and exclude this possibility. Thus we may assume that $p$ is odd.

We first establish an analogue of Theorem~2 in~\cite{LV2005} eliminating polytopes of high ranks.

\begin{lem}\label{highrank}
If $PGL_2(q)$ is the full automorphism group of a regular polytope $\cal P$, then the rank of $\cal P$ is at most 4.
\end{lem}
{\bf Proof:} 
Suppose $\cal P$ is a regular polytope of rank $n\geq 4$ with group $G :=\Gamma({\cal P}) = \langle \rho_{0},\ldots,\rho_{n-1} \rangle$. 
Then the subgroup $\langle\rho_{0},\rho_{2},\rho_{3},\ldots, \rho_{n-1}\rangle$ of $\Gamma({\cal P})$ is a subgroup of the centralizer $C_G(\rho_0)$ of $\rho_0$ in $G$.
By Proposition~\ref{subPGL}, there are two conjugacy classes of cyclic subgroups $Z_2$, of respective lengths $q(q-1)/2$ and $q(q+1)/2$. If $G$ acts by conjugation on the conjugacy class containing $\rho_0$, then clearly this action is transitive and the stabilizer of $\rho_0$ is its normalizer in $G$. But the normalizer of an involution is just its centralizer, and its index in $G$ is just the number of involutions in the conjugacy class. Therefore, the centralizer has index $q(q-1)/2$ or $q(q+1)/2$ in $G$ and thus it has order $2(q+1)$ or $2(q-1)$. Therefore, the rank is at most 4.
\hfill $\Box$
\smallskip

Alternatively we could have argued as follows. The subgroup 
$\langle\rho_{0},\rho_{1},\rho_{3},\rho_{4}\rangle$ of $\Gamma({\cal P})$ is of the form $D_{2k}\times D_{2l}$ for some $k,l$. Inspection of the list of subgroups of $PGL_2(q)$ given in Proposition~\ref{subPGL} then immediately shows that this cannot occur. While this argument is quick, the previous proof seems to have the advantage of potentially wider applicability to other types of groups.

\medskip
We next investigate the possibility for a regular polytope of rank $4$ to have a group of type $PGL_2(q)$. (As mentioned before, the case of rank~3 has already been settled in~\cite{SC94}.) Thus we assume that $(G, \{\rho_0,\ldots,\rho_{3}\})$ is a string C-group of type $\{t,l,s\}$ (that is, the orders of $\rho_0\rho_1$, $\rho_1\rho_2$ and $\rho_2\rho_3$ in $G$ are $t$, $l$ and $s$, respectively). Clearly, $t,l,s\geq 3$, since $G$ is not a direct product of two non-trivial groups. As before we let 
\[ G_i = \langle \rho_j \mid j \neq i\rangle \;\quad (\mathrm{for~} i=0,\ldots,3) . \]
If $\cal P$ is the regular $4$-polytope with group $G$, then $G_0$ must be the group of a vertex-figure and $G_3$ the group of a facet of $\cal P$. By Proposition~\ref{subPGL}, the only possible types of subgroups for $G_0$ and $G_3$ are $S_4$, $A_5$, the (Frobenius) groups of family (9) of Proposition~\ref{subPGL}, or $L_2(q')$ or $PGL_2(q')$ for some $q'$. Dihedral subgroups or $A_4$ cannot occur because they are not irreducible rank $3$ C-groups (see Lemma 2 in~\cite{LS2005a}). The following sequence of lemmas is aimed at eliminating the groups of family (9) as well as the subgroups $L_2(q')$ and $PGL_2(q')$ as possibilities. 

We require the following basic facts about the groups $PGL_{2}(q)$, which follow directly from similar properties for the projective special linear groups (see Lemma~3 in~\cite{LS2005a}); bear in mind here that $PGL_{2}(q)$ is a subgroup of $L_2(q^2)$.  The centre of a nonabelian subgroup of $G$ must have order at most 2. Moreover, if $H$ is a nonabelian subgroup of $G$ whose centre has size 2, then ($q$ is odd and) $H$ must be a dihedral group. 

\begin{lem}\label{lemma2}
The orders $t$ of $\rho_{0}\rho_{1}$ and $s$ of $\rho_{2}\rho_{3}$ must be odd.
\end{lem}
{\bf Proof:}
If $t$ is even, then the centre of the nonabelian subgroup 
\[ G_2=Z_{2}\times D_{2t}\cong (Z_{2}\times Z_{2})\times D_t \] 
is too large. Hence $t$, and similarly $s$, must be odd.
\hfill $\Box$
\smallskip

Our next two lemmas deal with subgroups of type $L_2(q')$ of $PGL_2(q)$. 

\begin{lem}\label{lemma2b}
Every subgroup $L_2(q')$ of $PGL_2(q)$ lies in the unique subgroup $L_2(q)$ of $PGL_2(q)$.
\end{lem}
{\bf Proof:}
This is implicit in the count of subgroups of type $L_{2}(q')$ for $PGL_2(q)$ and $L_{2}(q)$. In fact, there are $\frac{q(q^2-1)}{q'(q'^2-1)}$ subgroups $L_2(q')$ in $PGL_2(q)$ (see Proposition~\ref{subPGL}, item (7)) and the same number of subgroups $L_2(q')$ in $L_2(q)$ (see Proposition~\ref{sub}, item (10)).
\hfill $\Box$

\begin{lem}\label{lemma3}
Let $H$ and $K$ be two subgroups of type $L_2(q')$ in $PGL_2(q)$, with $q'^{\,m} = q$ for some positive integer $m$. 
Then $H\cap K$ cannot contain a dihedral group $D_{2k}$ with $k>2$ and $k$ a divisor of $\frac{(q'\pm 1)}{2}$.
\end{lem}
{\bf Proof:}
Since $H$ and $K$ are subgroups of $L_2(q)$, we may appeal to a similar such statement about subgroups of type $L_2(q')$ in $L_2(q)$ established in Lemma~6 of~\cite{LS2005a}. Lemma~6 itself is slightly weaker than the present lemma, but its proof actually establishes the stronger version presented here. For the convenience of the reader we repeat the argument here. (Note that the stronger version for $L_2(q)$ leads to some simplifications in the proofs of~\cite{LS2005a}.) 

Let $k>2$, and let $k$ be a divisor of $\frac{q'\pm 1}{2}$. From Proposition~\ref{sub} we
know that 
\begin{itemize}
\item in $L_2(q)$, there are $\frac{q(q^2-1)}{q'(q'^2-1)}$ subgroups isomorphic to 
$L_2(q')$ and $\frac{q(q^2-1)}{4k}$ subgroups isomorphic to $D_{2k}$;
\item in $L_2(q')$, there are $\frac{q'(q'^2-1)}{4k}$ subgroups isomorphic to
$D_{2k}$.
\end{itemize}
Let $n := \frac{(q'-1)}{2}$ if $k\mid \frac{(q'-1)}{2}$ and $n:=\frac{(q'+1)}{2}$ if $k\mid \frac{(q'+1)}{2}$.
By Proposition~\ref{sub}, there are $\frac{q(q^2-1)}{4n}$ subgroups $D_{2n}$ in $L_2(q)$. Each subgroup $D_{2n}$ contains $\frac{n}{k}$ subgroups $D_{2k}$.
Therefore, each subgroup $D_{2k}$ is contained in exactly one subgroup $D_{2n}$.
The same kind of arguments show that each $D_{2n}$ is contained in exactly one $L_2(q')$.
Therefore, each subgroup $D_{2k}$ of $L_2(q)$ (with $k>2$ and $k$ a divisor
of $\frac{(q'\pm 1)}{2}$) is contained in a subgroup $L_2(q')$ of $L_2(q)$, and the number of subgroups $L_2(q')$ containing a given subgroup $D_{2k}$ is precisely one. Now the lemma follows.
\hfill $\Box$
\smallskip

Observe that the same result may be proven in exactly the same way (with Proposition~\ref{sub} replaced by Proposition~\ref{subPGL}) if $G$ is $L_2(q)$. 

\begin{lem}\label{Frob}
The subgroups $G_0$ and $G_3$ of $G$ cannot be isomorphic to a group of family (9) of Proposition~\ref{subPGL}.
\end{lem}
{\bf Proof:} 
The involutions in a subgroup $E_{p^m}:Z_d$ of family (9) are exactly the elements of the form $\varphi\psi$ with $\varphi\in E_{p^m}$ and $\psi$ the involution in the cyclic factor $Z_d$; recall here that $E_{p^m}:Z_{d}$ is a subgroup of $E_q:Z_{q-1} \cong AGL_1(q)$. The product of two such involutions is necessarily of order $p$ or trivial. If $G_{3}$ (say) is isomorphic to $E_{p^m}:Z_d$, then this forces $\rho_0\rho_1$ to have order $p$. But then $\langle \rho_0,\rho_1,\rho_3\rangle$ is isomorphic to $D_{4p}$ which is impossible as we already mentioned earlier.
\hfill $\Box$

\begin{lem}\label{lemma4}
The subgroups $G_0$ and $G_3$ of $G$ cannot be isomorphic to $L_2(q')$, with
$q'^{\,m} = q$ for some positive integer $m$.
\end{lem}
{\bf Proof:} 
Suppose that $G_3$ (say) is isomorphic to $L_2(q')$. Then $G_3$ and its conjugate $G_3^{\rho_3}$ by $\rho_3$ are two subgroups of type $L_2(q')$ containing the dihedral group $G_{23} :=\langle\rho_0,\rho_1\rangle$. By Lemma~\ref{lemma3}, this dihedral group cannot be a group $D_{2k}$ with $k>2$ and $k$ a divisor of $\frac{(q'\pm 1)}{2}$.
So the only possibility left here is to have $G_{23} \cong D_{2p}$, with $p$ the prime divisor of $q$, which cannot occur.

Note here that $G_3$ and $G_3^{\rho_3}$ really are distinct subgroups. 
In fact, suppose that $G_3 = G_3^{\rho_3}$. 
Then, since $\rho_{0},\rho_{1},\rho_{2}\in G_3$, each generator $\rho_j$ of $G$ must normalize $G_3$ and hence $G_3$ must be normal in $G$. 
But $G_3$ is a subgroup of the simple group $L_{2}(q)$, so it must also be normal in $L_{2}(q)$ and hence coincide with $L_{2}(q)$.  
This forces the underlying polytope to have only two facets, which is impossible since $s>2$.
\hfill $\Box$
\smallskip

Next we eliminate the subgroups $PGL_2(q')$ as facet and vertex-figure groups.

\begin{lem}\label{lemma7}
The subgroups $G_0$ and $G_3$ of $G$ cannot be isomorphic to $PGL_2(q')$, with
$q'^{\,m}=q$ for some positive integer $m>1$.
\end{lem}
{\bf Proof:}
Suppose to the contrary that $G_3 \cong PGL_2(q')$ (say), with $q'^{\,m}=q$ and $m>1$. 
Then the two subgroups $G_3$ and $G_{3}^{\rho_3}$ are isomorphic to $PGL_2(q')$, and $D := G_{3}\cap G_{3}^{\rho_3}$ contains the dihedral group 
$D_{2t} := \langle\rho_{0},\rho_{1}\rangle$ of order $2t\,(\geq 6)$. The latter group is centralized by $\rho_3$. Note that $t\neq p$. 

Now, $G\cong PGL_2(q)$ can be embedded in a group $K \cong PGL_2(q^2)$. This group $K$ contains a subgroup $I$ of index 2 isomorphic to $L_2(q^2)$, and this subgroup $I$ itself contains $G$. Moreover, $G_3$ (resp. $G_3^{\rho_3}$) is contained in a subgroup $I_3$ (resp. $I_3^{\rho_3}$) isomorphic to $L_2(q'^2)$, and $I_3$ and $I_3^{\rho_3}$ are both contained in $I$. Since $I_3$ contains the element $\rho_0\rho_1$ of order $t$, we must have $t\mid \frac{(q'^{2} \pm 1)}{2}$, by Proposition~\ref{sub}, item (2), applied to $L_2(q'^2)$. 
Therefore, Lemma~\ref{lemma3} forces $I_3$ to be equal to $I_3^{\rho_3}$; in other words, $\rho_3$ must normalize $I_3$ in $K$.  Hence $\langle G_3, \rho_3\rangle = G \cong PGL_{2}(q)$ is a subgroup of the normalizer $N_K(I_3)$ of $I_3$ in $K$, and this normalizer is isomorphic to $PGL_2(q'^2)$. (Note for the latter that a supergroup of a subgroup $L_{2}(q'^2)$ of $PGL_{2}(q)$ is necessarily of the form $L_{2}(q'^{k})$ or $PGL_{2}(q'^{k})$ for some $k$; the original subgroup is normal in the supergroup only when $k=2$.) Therefore $m$ is at most $2$, and hence $m=2$ and $q=q'^{2}$. Now, $(\rho_2\rho_3)^s = 1$ and $s$ is odd, so $\rho_2$ and $\rho_3$ are conjugate in $G$.  But by Proposition~\ref{subPGL}, item (1), applied to $G\cong PGL_2(q)$, since $\rho_2$ is contained in $I_3\cong L_2(q)$, so is $\rho_3$. It follows that $\langle \rho_0,\rho_1,\rho_2,\rho_3\rangle = L_2(q) < G$, a contradiction. 
\hfill $\Box$
\medskip

Finally, then, we obtain our main theorem which was first conjectured to be true based on computations for the atlas of polytopes~\cite{LV2005}.

\begin{thm}\label{ourthm}
If $PGL_2(q)$ is the full automorphism group of a regular polytope $\cal P$ of rank $n\geq 4$, then $n=4$, $q=5$, and $\cal P$ is the $4$-simplex $\{3,3,3\}$. Rephrased in terms of C-groups, if $(G, \{\rho_0,\ldots, \rho_{n-1}\})$ is a string C-group of rank $n\geq 4$ isomorphic to $PGL_2(q)$, then $n=4$, $q=5$, and $G\cong S_{5}$ (occurring in its natural representation as the group of the $4$-simplex). 
\end{thm}
{\bf Proof:}
Lemmas~\ref{Frob}, ~\ref{lemma4} and~\ref{lemma7} (as well as Lemma~\ref{highrank}) reduce the possible types of subgroups of $G$ that can occur as $G_0$ and $G_3$ to only two kinds in each case, namely $S_4$ and $A_5$. We show that this leaves only one possibility.

It is well-known (and straightforward to check) that the only rank 3 polytopes with group $S_4$ are the tetrahedron $\{3,3\} \,(=\{3,3\}_4)$, the hemi-octahedron $\{3,4\}_3$, and the hemi-cube $\{4,3\}_3$; and those with group $A_5$ are the hemi-icosahedron $\{3,5\}_5$, the hemi-dodecahedron $\{5,3\}_5$, and the great dodecahedron $\{5,5\}_3 \,(=\{5,\frac{5}{2}\})$. By \cite{LS2005a}, Table~1, we readily see that, when their automorphism groups are taken in pairs to form the vertex-figure group $G_0$ and facet group $G_3$ of the group $G$ of a regular rank $4$ polytope, then $G$ can only be isomorphic to a group $PGL_{2}(q)$ if both $G_0$ and $G_3$ are groups $S_4$ in its natural representation as the group of a tetrahedron $\{3,3\}$, that is, if $G$ itself is the group $PGL_2(q) \cong S_5$ in its natural representation as the group of the $4$-simplex $\{3,3,3\}$.
\hfill $\Box$
\bigskip

Our main result can be rephrased by saying that the $4$-simplex is the only abstract polytope of rank $4$ (or higher) on which a group of type $PGL_{2}(q)$ can admit a faithful transitive action on the flags. It is quite remarkable that no similar result can hold for actions with two flag orbits. In fact, there are many chiral $4$-polytopes whose automorphism group is of type $PGL_{2}(q)$. (Recall here that a polytope is chiral if it has two orbits on the flags such that any two adjacent flags are in distinct orbits.) For example, if $p$ is a prime with $p \equiv 5\, (8)$ and $b$, $c$ are positive integers with $p=b^{2}+c^{2}$, then there exists a chiral $4$-polytope $\cal P$ of type $\{4,4,3\}$ with toroidal facets $\{4,4\}_{(b,c)}$ and with group isomorphic to $PGL_{2}(p)$ (see \cite[p.239,240]{SW1994}). For $p=5$ and $p=13$ these are the universal chiral polytopes $\{\{4,4\}_{(2,1)},\{4,3\}\}$ and $\{\{4,4\}_{(3,2)},\{4,3\}\}$ with groups $PGL_{2}(5)$ and $PGL_{2}(13)$, respectively (see also \cite{colw}). Similar examples also exist for other Schl\"afli symbols.
\medskip

As observed in \cite{SW1994} (see also Section 6 of~\cite{LS2005a}), if $p$ is a prime with $p\equiv 3(4)$, then there are regular $4$-polytopes whose rotation (even) subgroup is isomorphic to $L_{2}(p^2)$ and has index $2$ in the full group. Then, by Theorem~\ref{ourthm}, these polytopes must have a group different from $PGL_{2}(p^2)$. 
In fact, given an odd square prime power $q$, the outer automorphism group of $L_2(q)$ contains a Klein group $D_4 = \{1, a, b, c\}$, such that the extension $L_2(q)\langle a\rangle$ of $L_{2}(q)$ is isomorphic to $PGL_2(q)$ and such that $b$ is induced by the field automorphism of order 2. Then the extension $L_2(q)\langle c\rangle$ is sharply 3-transitive on the points of the projective line, as is $PGL_2(q)$. There is a third subgroup of index 2 in the automorphism group $\rm{Aut}(L_2(q))$, namely the extension $L_2(q)\langle b\rangle$ usually denoted by $P\Sigma L_2(q)$. This group is not sharply 3-transitive on the points of the projective line.

For $q=9$, $P\Sigma L_2(9) \cong S_6$, which is known to occur as group of regular polytopes of ranks $4$ and $5$ (see~\cite{LV2005}). The group $L_2(9)\langle c\rangle$ is isomorphic to the Mathieu group $M_{10}$ and is known to be a group that cannot be generated by involutions. Hence it is not the group of a regular polytope.

For $q=25$, a computer search produces, up to isomorphism, 17 rank four polytopes for $P\Sigma L_2(25)$, with Schl\"afli symbols 
$\{3,3,6\}$, $\{3,4,5\}$, $\{3,5,3\}$, $\{3,6,3\}$, $\{3,6,4\}$, $\{3,6,5\}$, $\{3,6,6\}$, $\{4,3,5\}$, $\{4,3,6\}$,
$\{4,5,4\}$, $\{4,5,5\}$, $\{4,6,4\}$, $\{5,3,6\}$, $\{5,4,5\}$, $\{5,5,6\}$, $\{5,6,5\}$ and $\{6,4,6\}$.
Again, the group $L_2(25)\langle c\rangle$ has no polytope.

\begin{conj}
Let $q$ be an odd square prime power, and let $L_2(q)\langle c\rangle$ and $L_2(q)\langle b\rangle = P\Sigma L_2(q)$ be as defined above. Then, \\
(a) $L_2(q)\langle c\rangle$ is not the group of a regular polytope;\\
(b) $L_2(q)\langle b\rangle$ is the group of a regular polytope of rank $4$.
\end{conj}

We remark that groups of type $L_{2}(q)$ or $PGL_{2}(q)$ also occur quite frequently as quotients of certain hyperbolic Coxeter groups (or their rotation subgroups) with non-string diagrams (see~\cite{CL2007},~\cite{JL2008}).

{\bf Acknowledgements:} The authors thank the referees and Daniel Pellicer for very useful comments on the preliminary version of this paper.

\bibliographystyle{plain}

\end{document}